\documentclass{groups13}
\usepackage{verbatim}
\newtheorem{conj}[thm]{Conjecture}
\begin{document}


%
\maketitle{RECENT WORK ON BEAUVILLE SURFACES, STRUCTURES AND GROUPS}{BEN FAIRBAIRN}{Department of Economics, Mathematics and Statistics, Birkbeck, University of London, Malet Street, London, WC1E 7HX\newline
Email: b.fairbairn@bbk.ac.uk}

\begin{abstract}
Beauville surfaces are a class of complex surfaces defined by letting a finite group
$G$ act on a product of Riemann surfaces. These surfaces possess many attractive
geometric properties several of which are dictated by properties of the group $G$.
In this survey we discuss the groups that may be used in this way. \emph{En route} we
discuss several open problems, questions and conjectures.
\end{abstract}

\section{Introduction}
Roughly speaking (precise definitions will be given in the next section), a Beauville
surface is a complex surface $\mathcal{S}$ defined by taking a pair of complex curves, i.e.
Riemann surfaces, $\mathcal{C}_1$ and $\mathcal{C}_2$ and letting a finite group $G$ act freely on their product to
define $\mathcal{S}$ as a quotient $(\mathcal{C}_1\times\mathcal{C}_2)/G$. These surfaces have a wide variety of attractive
geometric properties: they are surfaces of general type; their automorphism groups
\cite{41} and fundamental groups \cite{16} are relatively easy to compute (being closely
related to $G$ - see Section \ref{F1} and \ref{F2}); these surfaces are rigid surfaces in the sense of admitting no nontrivial
deformations \cite{9} and thus correspond to isolated points in the moduli space of
surfaces of general type \cite{30}.

Much of this good behaviour stems from the fact that the surface $(\mathcal{C}_1\times\mathcal{C}_2)/G$
is uniquely determined by a particular pair of generating sets of $G$ known as a `Beauville
structure'. This converts the study of Beauville surfaces to the study of groups
with Beauville structures, i.e. Beauville groups.

Beauville surfaces were first defined by Catanese in \cite{16} as a generalisation of an
earlier example of Beauville \cite[Exercise X.13(4)]{13} (native English speakers may find the English translation \cite{14} somewhat easier to read and get hold of) in which $\mathcal{C}=\mathcal{C}'$ and the curves
are both the Fermat curve defined by the equation $X^5 + Y^5 + Z^5 = 0$ being acted
on by the group $(\mathbb{Z}/5\mathbb{Z})\times(\mathbb{Z}/5\mathbb{Z})$ (this choice of group may seem somewhat odd
at first, but the reason will become clear later). Bauer, Catanese and Grunewald
went on to use these surfaces to construct examples of smooth regular surfaces with
vanishing geometric genus \cite{10}. Early motivation came from the consideration of
the `Friedman-Morgan speculation' -- a technical conjecture concerning when two
algebraic surfaces are diffeomorphic which Beauville surfaces provide counterexamples to. More recently, they have been used to construct interesting orbits of the
absolute Galois group Gal($\overline{\mathbb{Q}}/\mathbb{Q}$) (connections with Gothendeick's theory of dessins
d'enfant make it possible for this group to act on the set of all Beauville surfaces).
We will discuss this in slightly more detail in Section 7.6. Furthermore, Beauville's
original example has also recently been used by Galkin and Shinder in \cite{27} to construct
examples of exceptional collections of line bundles.

Like any survey article, the topics discussed here reflect the research interests
of the author. Slightly older surveys discussing related geometric and topological
matters are given by Bauer, Catanese and Pignatelli in \cite{11, 12}. Other notable
works in the area include \cite{7,42,CombSurv,Wo}.

We remark that throughout we shall use the standard `Atlas' natation for finite groups and related concepts as described in \cite{ATLAS}, excepting that we will occasionally deviate to minimise confusion with similar notation for geometric concepts.

In Section 2 we will introduce the preliminary definitions before proceeding in
Section 3 to discuss the case of the finite simple groups. We then go on in Section 4
to discuss the abelian and nilpotent groups. Next, we focus our attention on special
types of Beauville structures when we discuss strongly real Beauville structures
in Section 5 and mixed Beauville structures in Section 6. Finally, we discuss a
miscellany of related but less well studied topics in Section 7.

\section{Preliminaries}
\begin{de}
A surface $\mathcal{S}$ is a \textbf{Beauville surface of unmixed type} if
\begin{itemize}
\item the surface $\mathcal{S}$ is isogenous to a higher product, that is, $\mathcal{S}\cong(\mathcal{C}_1\times\mathcal{C}_2)/G$ where
$\mathcal{C}_1$ and $\mathcal{C}_2$ are algebraic curves of genus at least 2 and $G$ is a finite group acting
faithfully on $\mathcal{C}_1$ and $\mathcal{C}_2$ by holomorphic transformations in such a way that it
acts freely on the product $\mathcal{C}_1\times\mathcal{C}_2$, and
\item each $\mathcal{C}_i/G$ is isomorphic to the projective line $\mathbb{P}_1(\mathbb{C})$ and the covering map
$\mathcal{C}_i\rightarrow\mathcal{C}_i/G$ is ramified over three points.
\end{itemize}
\end{de}

There also exists a concept of Beauville surfaces of mixed type but we shall
postpone our discussion of these until Section 6.
In the first of the above conditions the genus of the curves in question needs to
be at least 2. It was later proved by Fuertes, Gonz\'{a}lez-Diez and Jaikin-Zapirain in \cite{25} that
in fact we can take the genus as being at least 6.
The second of the above conditions implies that each $\mathcal{C}_i$ carries a regular dessin
in the sense of Grothendieck's theory of \emph{dessins d'enfants} (childeren's drawings) \cite{36}. Furthermore, by Bely\u{\i}'s
Theorem \cite{15} this ensures that $\mathcal{S}$ is defined over an algebraic number field in the sense that when we
view each Riemann surface as being the zeros of some
polynomial we find that the coeffcients of that polynomial belong to some number field.
Equivalently they admit an orientably regular hypermap \cite{43}, with $G$ acting as
the orientation-preserving automorphism group. A modern account of dessins d'enfants and proofs of Bely\u{\i}'s theorem may be found in the recent book of Girondo and Gonz\'{a}lez-Diez \cite{GGDbook}.

This can also be described instead in terms of uniformisation and the language
of Fuchsian groups \cite{31, 48}.

What makes this class of surfaces so good to work with is the fact that all of the
above definition can be `internalised' into the group. It turns out that a group $G$
can be used to define a Beauville surface if and only if it has a certain pair of generating
sets known as a Beauville structure.

\begin{de}\label{MainDef} Let $G$ be a finite group. Let $x,y\in G$ and let
$$\Sigma(x, y) :=\bigcup_{i=1}^{|G|}\bigcup_{g\in G}\{(x^i)^g,(y^i)^g,((xy)^i)^g\}.$$

An \textbf{unmixed Beauville structure} for the group $G$ is a set of pairs of elements $\{\{x_1, y_1\},\{x_2,y_2\}\}\subset G\times G$
with the property that $\langle x_1, y_1\rangle = \langle x_2, y_2\rangle=G$ such that
$$\Sigma(x_1, y_1)\cap \Sigma(x_2, y_2)=\{e\}.$$
If $G$ has a Beauville structure we say that $G$ is a \textbf{Beauville group}. Furthermore
we say that the structure has \textbf{type}
\[
((o(x_1), o(y_1), o(x_1y_1)),(o(x_2), o(y_2), o(x_2y_2))).
\]
\end{de}

Traditionally, authors have defined the above structure in terms of so-called `spherical systems of generators of length 3', meaning $\{x, y, z\}\subset G$ with $xyz = e$, but we
omit $z = (xy)^{-1}$ from our notation in this survey. (The reader is warned that this terminology is a little misleading since the underlying geometry of Beauville surfaces is hyperbolic thanks to the below constraint on the orders of the elements.) Furthermore, many earlier
papers on Beauville structures add the condition that for $i = 1, 2$ we have that
\[
\displaystyle \frac{1}{o(x_i)}+\frac{1}{o(y_i)}+\frac{1}{o(x_iy_i)}<1,
\]
but this condition was subsequently found to be
unnecessary following Bauer, Catanese and Grunewald's investigation of the wall-paper groups in \cite{8}. A triple of elements and their orders satisfying this condition
are said to be hyperbolic. Geometrically, the type gives us considerable amounts of
geometric information about the surface: the Riemann-Hurwitz formula
\[
\displaystyle g(\mathcal{C}_i)=1+\frac{|G|}{2}\bigg(1-\frac{1}{o(x_i)}-\frac{1}{o(y_i)}-\frac{1}{o(x_iy_i)}\bigg)
\]
tells us the genus of each of the curves used to define the surface $\mathcal{S}$ and by a theorem
of Zeuthen-Segre this in turn gives us the Euler number of the surface $\mathcal{S}$ since
\[
\displaystyle e(\mathcal{S})=4\frac{(g(\mathcal{C}_1)-1)(g(\mathcal{C}_2)-1)}{|G|}
\]
which in turn gives us the holomorphic Euler-Poincar\'{e} characteristic of $\mathcal{S}$, namely
$4\chi(\mathcal{S})=e(\mathcal{S})$ (see \cite[Theorem 3.4]{16}).

Furthermore, if a group can be generated by a pair of elements of orders $a$ and
$b$ whose product has order $c$ then $G$ is a homomorphic image of the triangle group
$$T_{a,b,c} = \langle x, y, z | x^a = y^b = z^c = xyz = 1\rangle.$$
Homomorphic images of the triangle group $T_{2,3,7}$ are known as Hurwitz groups.
In several places in the literature, knowing that a particular group is a Hurwitz
group has proved useful for deciding its status as a Beauville group. For a discussion of known results on Hurwitz groups see the excellent surveys of Conder \cite{Conder1,18}.

\section{Finite Simple Groups}

A necessary condition for a group to be a Beauville group is that it is 2-generated.
In \cite{1,47} it is proved that all non-abelian finite simple groups are 2-generated. For a
long time it was conjectured that every non-abelian finite simple group, aside from
the alternating group A$_5$, is a Beauville group \cite[Conjecture 7.17]{9}, providing a
rich source of examples. Various authors proved special cases of this \cite{9,24,26}.
The full result comes from the following Theorem which is proved by the author,
Magaard and Parker in \cite{21, 22}.
\begin{thm}\label{SimpThm} With the exceptions of $SL_2(5)$ and $PSL_2(5)(\cong A_5\cong
SL_2(4))$, every
finite quasisimple group is a Beauville group.
\end{thm}
Similar results were proved at around the same time by Garion, Larsen and Lubotzky in \cite{29} (using probabilistic results concerning triangle groups from the PhD thesis of Marion \cite{44}) and
by Guralnick and Malle in \cite{37} using the theory of linear algebraic groups. Since
the overriding ideas behind the proofs given in \cite{21, 29, 37} are in many ways quite
general we sketch these ideas in the hope that they may be useful in proving other
conjectures that appear later in this survey.

First note that the alternating groups can be dealt with using classical permutation group theory. Furthermore, the low rank groups of Lie type may be dealt with
using explicit matrix calculations (see for instance the work of Fuertes and Jones
in \cite{26} concerning the groups $PSL_2(q)$, $^2B_2(2^{2n+1})$ and $^2G_2(3^{2n+1})$.) The sporadic
simple groups are easily dealt with on a case by case basis with structure constant
calculations being useful for the larger groups. The real difficulty lies with the
groups of Lie type of unbounded rank.

Let $G$ be a finite simple group of Lie type of characteristic $p$. To ensure that
we can choose elements of the group $G$ whose product behaves as we require we use a theorem of
Gow \cite{35} (a slight generalisation of this result to quasisimple groups is given in \cite[Theorem 2.6]{21}). An element of $G$ is said to be `semisimple' if its order is coprime to $p$
and is said to be `regular semisimple' if its centralizer in $G$ has order coprime to $p$.

\begin{thm} Let $G$ be a finite simple group of Lie type of characteristic $p$ and
let $s\in G$ be a semisimple element. Let $R_1, R_2\subset G$ be conjugacy classes of regular
semisimple elements of $G$. Then there exist elements $x\in R_1$ and $y\in R_2$ such that
$s = xy$.
\end{thm}

To ensure that the conjugacy part of the definition of a Beauville structure is
satisfied we aim to choose $x_1, x_2, y_1, y_2\in G$ such that $o(x_1)o(y_1)o(x_1y_1)$ is coprime
to $o(x_2)o(y_2)o(x_2y_2)$. This is made possible by a classical theorem of Zsigmondy
\cite{50} (or rather Bang \cite{2} in the case $p=2$.) Whilst \cite{2} and \cite{50} are over a century old and therefore difficult to read and get hold of, a more recent account of a proof is given by L\"{u}nburg in \cite{Lub}.

\begin{thm}\label{Zig} For any positive integers $a$ and $n$ there exists a prime that divides
$a^n-1$ but not $a^k-1$ for any $k < n$ with the following exceptions:
\begin{itemize}
\item $a = 2$ and $n = 6$; and
\item $a + 1$ is a power of $2$ and $n = 2$.
\end{itemize}
\end{thm}

The real significance of the above results stems from the fact that most groups of
Lie type have an order that is a product of numbers of the form $p^k-1$ and so the
above result guarantees the existence of a rich supply of distinct primes that can
be taken as being the orders of the elements of our Beauville structure.

It remains to decide if a given triple will generate the group. Since our elements
have orders given by Theorem \ref{Zig} we can use a theorem of Guralnick, Pentilla,
Praeger and Saxl \cite{38} concerning subgroups of the general linear group $GL_n(p^a)$ containing elements of
these orders and closely related results of Niemeyer and Praeger \cite{45} for the other
classical groups to show that no proper subgroups contain our elements. It follows that our
chosen elements will generate the group.

\section{Abelian and Nilpotent Groups}

The abelian Beauville groups were essentially classified by Catanese in \cite[page 24.]{16} and the full argument is given explicitly in \cite[Theorem 3.4]{8} where the following
is proved.

\begin{thm}\label{AbThm} Let $G$ be an abelian group. Then $G$ is a Beauville group if, and only if, $G = (\mathbb{Z}/n\mathbb{Z})\times(\mathbb{Z}/n\mathbb{Z})$ where $n>1$ is
coprime to $6$.
\end{thm}

This explains why Beauville's original example used the group $(\mathbb{Z}/5\mathbb{Z})\times(\mathbb{Z}/5\mathbb{Z})$ - it is
the smallest abelian Beauville group.

Theorem \ref{AbThm} has been put to great use by Gonz\'{a}lez-Diez, Jones and Torres-Teigell in \cite{33} where several structural results concerning the surfaces defined by
abelian Beauville groups are proved. For each abelian Beauville group they describe
all the surfaces arising from that group, enumerate them up to isomorphism and
impose constraints on their automorphism groups. As a consequence they show
that all such surfaces are defined over $\mathbb{Q}$.

After the abelian groups, the next most natural class of finite groups to consider
are the nilpotent groups. In \cite[Lemma 1.3]{3} Barker, Boston and the author note
the following easy Lemma.

\begin{lemma}
Let $G$ and $G'$ be Beauville groups and let $\{\{x_1,y_1,\},\{x_2,y_2\}\}$ and $\{\{x'_1,y'_1,\},\{x'_2,y'_2\}\}$ be their respective Beauville structures. Suppose that
$$gcd(o(x_i),o(x_i'))=gcd(o(y_i),o(y_i'))=1$$
for $i=1,2$. Then $\{\{(x_1,x'_1,),(y_1,y'_1)\},\{(x_2,x'_2),(y_2,y'_2)\}\}$ is a Beauville structure for the group $G\times G'$.
\end{lemma}

Recall that a finite group is nilpotent if, and only if, it isomorphic to the direct product of its Sylow subgroups.
It thus follows that this lemma, and its easy to prove converse, reduces the study of nilpotent
Beauville groups to that of Beauville $p$-groups. Note that Theorem \ref{AbThm} gives us
infinitely many examples of Beauville $p$-groups for every prime $p > 3$ - simply let
$n$ be any power of $p$. Early examples of Beauville 2-groups and 3-groups
were constructed by Fuertes, Gonz\'{a}lez-Diez and Jaikin-Zapirain in \cite{25} where a
Beauville group of order $2^{12}$ and another of order $3^{12}$ were constructed. Even earlier than this, two Beauville 2-groups of order $2^8$ arose as part of a classification due to Bauer, Catanese and Grunewald in \cite{10} of certain classes of surfaces of general type.

More recently, in \cite{3} Barker, Boston and the author classified the Beauville
$p$-groups of order at most $p^4$ and made substantial progress on the cases of groups
of order $p^5$ and $p^6$. In particular, the number of Beauville $p$-groups of order
$p^4$ is two for every $p > 3$ and zero otherwise, but for $p^5$ we have the following.

\begin{conj}
For all $p\geq5$, the number of Beauville $p$-groups of order $p^5$ is
given by $p + 10$.
\end{conj}

In \cite[Theorem 1.4]{3} we prove that there are at least $p + 8$ Beauville groups of
order $p^5$. Furthermore, the above conjecture has been verified computationally for
all primes $p$ such that $5\leq p\leq 19$. Perhaps more interestingly, other results proved
in \cite{3} verify that the proportion of 2-generated $p$-groups of order $p^5$ that are Beauville tends to
1 as $p$ tends to infinity, however this fails to to be true for $p$-groups of order $p^6$.

\begin{qu} If $n > 6$ what is the behaviour, as $p$ tends to infinity, of the proportion of 2-generated $p$-groups that are Beauville?
\end{qu}

Another consequence of this work was determining the smallest Beauville $p$-group for all primes. In the below presentations, if no relationship between two generators is specified by a relation or relator then it should be assumed that the two generators commute.

\begin{thm} The smallest Beauville $p$-groups are as follows.
\begin{itemize}
\item For $p = 2$ the group
$$\langle x, y \,|\, x^4, y^4, [x,y^2]^22, [x,y^3]^2, [x^2, y^3]\rangle$$
of order $2^7$.
\item For $p = 3$ the group
$$\langle x, y, z, w, t \,|\, x^3, y^3, z^3, w^3, t^3, y^x = yz, z^x = zw, z^y = zt\rangle$$
of order $3^5$.
\item For $p\geq5$ the group
$$\langle x, y, z \,|\, x^5, y^5, z^5, [x, y] = z\rangle$$
of order $p^3$.
\end{itemize}
\end{thm}

Further examples are given by the following unpublished constructions due to
Jones and Wolfart.

\begin{thm} Let $G$ be a finite group of exponent $n = p^e > 1$ for some prime
$p\geq5$, such that the abelianisation $G/G'$ of $G$ is isomorphic to $\mathbb{Z}_n\times\mathbb{Z}_n$. Then $G$
has a Beauville structure.
\end{thm}

\begin{co} Let $G$ be a $2$-generated finite group of exponent $p$ for some prime
$p\geq5$. Then $G$ has a Beauville structure.
\end{co}

As noted earlier Beauville $p$-groups for $p>3$ are in bountiful supply. Several
examples of Beauville 2-groups and 3-groups are constructed by Barker, Boston,
Peyerimhoff and Vdovina in \cite{5, 6} by considering sections of groups defined using
projective planes. More recently, in \cite{4} Barker, Boston, Peyerimhoff and Vdovina
using similar ideas constructed the first infinite family of Beauville 2-groups. At
the time of writing, as far as the author is aware, only finitely many Beauville
3-groups are known leading to the following natural problem.

\begin{prob}
Construct infinitely many Beauville 3-groups.
\end{prob}

More recently in \cite{SV} Stix and Vdovina give a construction of Beauville $p$-groups that provides infinitely many examples for every $p\geq5$. More specifically they prove the following.
\begin{thm}
Let $p$ be a prime, $n,m\in\mathbb{N}$ and $\lambda\in(\mathbb{Z}/p^m\mathbb{Z})^{\times}$ with $\lambda^{p^n}\cong1$ mod $p^m$. The semidirect product $\mathbb{Z}/p^m\mathbb{Z}:\mathbb{Z}/p^n\mathbb{Z}$ with action $\mathbb{Z}/p^n\mathbb{Z}\rightarrow\mbox{Aut}(\mathbb{Z}/p^m\mathbb{Z})$ sending $1\mapsto\lambda$ admits an unmixed Beauville structure if and only if $p\geq5$ and $n=m$.
\end{thm}
They go on to prove related results using the theory of pro-$p$ groups.

We conclude this section with the following remarks. Nigel Boston has recently
undertaken some substantial and as yet unpublished computations regarding the
relationship between $p$-groups' status as Beauville groups and their position on the so-called
`O'Brien Trees' \cite{46}. Whilst little global pattern appears to exist in general, there
does appear to be some mysterious relationship with an invariant known as the
`nuclear rank' of the group - see \cite{Boston}. Since defining this concept is somewhat involved we
shall say no more about this here.

\section{Strongly Real Beauville Groups}

Given any complex surface $\mathcal{S}$ it is natural to consider the complex conjugate surface
$\overline{\mathcal{S}}$. In particular it is natural to ask if the surfaces are biholomorphic.

\begin{de} Let $\mathcal{S}$ be a complex surface. We say that $\mathcal{S}$ is \textbf{real} if there exists
a biholomorphism $\sigma:\mathcal{S}\rightarrow\overline{\mathcal{S}}$ such that $\sigma^2$ is the identity map.
\end{de}

As noted earlier this geometric condition can be translated into algebraic terms.

\begin{de}
 Let $G$ be a Beauville group and let $X =\{\{x_1, y_1\},\{x_2, y_2\}\}$ be a
Beauville structure for $G$. We say that $G$ and $X$ are \textbf{strongly real} if there exists an
automorphism $\phi\in\mbox{Aut}(G)$ and elements $g_i\in G$ for $i = 1, 2$
such that
$$g_1\phi(x_i)g_1^{-1}=x_i^{-1}\mbox{ and }g_2\phi(y_i)g_2^{-1}=y_i^{-1}$$
for $i=1,2$.
\end{de}

It is often, but not always, convenient to take $g_1=g_2$.

Our first examples come immediately from Theorem \ref{AbThm} since for any abelian
group the function $x\mapsto-x$ is an automorphism.

\begin{co}\label{AbCor} Every Beauville structure of an abelian Beauville group is strongly
real.
\end{co}

A little more generally, when it comes to strongly real Beauville $p$-groups the
examples given by Theorem \ref{AbThm} are, as far as the author is aware, the only known examples. Furthermore,
the Beauville 2-groups constructed by Barker, Boston, Peyerimhoff and Vdovina in \cite{4} are explicitly shown to not be strongly real. However, a combination of Corollary \ref{AbCor} and the fact that $p$-groups in general tend to have large automorphism groups \cite{Boston2,BostonBushHajir}
it seems likely that most Beauville $p$-groups are in fact strongly real Beauville
groups. This makes the following problem
particularly pressing.

\begin{prob}
Construct examples of strongly real Beauville $p$-groups.
\end{prob}

In \cite{19} the following conjecture, a refinement of an earlier conjecture of Bauer,
Catanese and Grunewald \cite[Section 5.4]{8}, is made.

\begin{conj} All non-abelian finite simple groups apart from A$_5$, M$_{11}$ and M$_{23}$
are strongly real Beauville groups.
\end{conj}

Only a few cases of this conjecture are known.
\begin{itemize}
\item In \cite{24} Fuertes and Gonz\'{a}lez-Diez showed that the alternating groups $A_n$
($n\geq 7$) and the symmetric groups $S_n$ ($n\geq 5$) are strongly real Beauville
groups by explicitly writing down permutations for their generators and the
automorphisms used and applying some of the classical theory of permutation
groups to show that their elements had the properties they claimed. It was subsequently found that the group A$_6$ is also strongly real.
\item In \cite{26} Fuertes and Jones proved that the simple groups $PSL_2(q)$ for prime
powers $q > 5$ and the quasisimple groups $SL_2(q)$ for prime powers $q > 5$
are strongly real Beauville groups. As with the alternating and symmetric
groups, these results are proved by writing down explicit generators, this time
combined with a celebrated theorem usually (but historically inaccurately)
attributed to Dickson for the maximal subgroups of $PSL_2(q)$. (For a full statement of this result and related theorems as well a detailed historical account of the maximal subgroups of low dimensional classical groups see the excellent survey of King in \cite{King}.) General lemmas
for lifting structures from a group to its covering groups are also used.
\item  In \cite{20} the author determined which of the sporadic simple groups are strongly
real Beauville groups, including the `27$^{th}$ sporadic simple group', the Tits
group $^2F_4(2)'$. Only the Mathieu groups M$_{11}$ and M$_{23}$ are not strongly real.
For all of the other sporadic groups smaller than the Baby Monster group $\mathbb{B}$
explicit words in the `standard generators' \cite{49} for a strongly real Beauville
structure were given. For the Baby Monster group $\mathbb{B}$ and Monster group $\mathbb{M}$
character theoretic methods were used.
\item In \cite{19} the author also verified this conjecture for the Suzuki groups $^2B_2(2^{2n+1})$.
Again, this was achieved by writing down explicit elements of the group which
using the list of maximal subgroups of the Suzuki group are shown to generate.
\item  In \cite{19} the author extended earlier computations of Bauer, Catanese and
Grunewald, verifying this conjecture for all non-abelian finite simple groups
of order at most $100\,000\,000$.
\end{itemize}

We remark that several of the groups mentioned in the above bullet points are
not simple. More generally we ask the following.

\begin{qu}
Which groups are strongly real Beauville groups?
\end{qu}

Finally, we remark that in \cite{19} the author constructs many further examples of strongly real Beauville groups. This includes the characteristically simple groups $A_n^k$ for moderate values of $k$ and sufficiently large values of $n$, the groups $S_n\times S_n$ for $n\geq5$ and the almost simple sporadic groups. This last calculation combined with the earlier remarks on the symmetric group lead to the following conjecture.

\begin{co}
A split extension of a simple group is a Beauville group if, and
only if, it is a strongly real Beauville group.
\end{co}

\section{The Mixed Case}

When we defined Beauville surfaces and groups we considered the action of a group
$G$ on the product of two curves $\mathcal{C}_1\times\mathcal{C}_2$. In an unmixed structure this action comes
solely from the action of $G$ on each curve individually, however there is nothing to
stop us considering an action on the product that interchanges the two curves and
it is precisely this situation that we discuss in this section. Recall from Definition
\ref{MainDef} that given $x, y\in G$ we write
$$\Sigma(x, y) :=\bigcup_{i=1}^{|G|}\bigcup_{g\in G}\{(x^i)^g,(y^i)^g,((xy)^i)^g\}.$$

\begin{de} Let $G$ be a finite group. A \textbf{mixed Beauville structure} for $G$ is
a quadruple $(G^0, g, h, k)$ where $G^0$ is an index $2$ subgroup and $g, h, k\in G$ are such
that
\begin{itemize}
\item$\langle g, h\rangle = G^0$;
\item$k\not\in G^0$;
\item for every $\gamma\in G^0$ we have that $(k\gamma)^2\not\in\Sigma(g,h)$ and
\item $\Sigma(g,h)\cap\Sigma(g^k,h^k) =\{e\}$
\end{itemize}
A Beauville surface defined by a mixed Beauville structure is called a \textbf{mixed
Beauville surface} and group possessing a mixed Beauville structure is called a \textbf{mixed Beauville group}.
\end{de}

In terms of the curves defining the surface, the group $G^0$ is the stabiliser of the
curves with the elements of $G\setminus G^0$ interchanging the two terms of $\mathcal{C}_1\times\mathcal{C}_2$. Moreover
it is only possible for a Beauville surface $(\mathcal{C}_1\times\mathcal{C}_2)/G$ to come from a mixed Beauville
structure if $\mathcal{C}_1\cong\mathcal{C}_2$. The above conditions also ensure that $\{\{g,h\},\{g^k,h^k\}\}\subset G^0\times G^0$ is a
Beauville structure for $G^0$.

In general, mixed Beauville structures are much harder to construct than their
unmixed counterparts. The following lemma of Fuertes and Gonz\'{a}lez-Diez imposes
a strong condition on a group with a mixed Beauville structure \cite[Lemma 5]{24}.

\begin{lemma} Let $(\mathcal{C}_1\times\mathcal{C}_2)/G$ be a mixed Beauville surface and let $G^0$ be the subgroup
of $G$ consisting of the elements which do not interchange the two curves. Then the
order of any element in $G\setminus G^0$ is divisible by $4$.
\end{lemma}

Clearly no simple group can have a mixed Beauville structure since it is necessary
to have a subgroup of index $2$ and the cyclic group of order $2$ is not a Beauville
group, however that does not preclude the possibility of almost simple groups having
mixed Beauville structures. The above lemma was originally used to show that
no symmetric group has a mixed Beauville structure. In \cite{20} the author used the
above to show that no almost simple sporadic group has a mixed Beauville structure
(though the almost simple Tits group $^2F_4(2)$ is not excluded by the above lemma) and in general most
almost simple groups are ruled out by it (though as the groups $P\Sigma L_2(p^2)$ show there are infinitely many exceptions to this). A further restriction comes from \cite[Theorem 4.3]{8} where Bauer, Catanese and Grunewald prove that $G^0$ must be non-abelian. Various geometric
constraints are proved by Torres-Teigell in his PhD thesis \cite{48}. Most notably the genus of a mixed Beauville surface is odd and at least 17. Furthermore, this bound is sharp. This naturally leads to the following problem.

\begin{prob} Find mixed Beauville structures.
\end{prob}

The earliest examples of groups that do possess mixed Beauville structures were
given by Bauer, Catanese and Grunewald in \cite{8}. Their general construction is of
the form $(H\times H):(\mathbb{Z}/4\mathbb{Z})$, the generator of the group $\mathbb{Z}/4\mathbb{Z}$ acting on the direct
product by interchanging its two factors and $G^0=H\times H\times \mathbb{Z}/2\mathbb{Z}$.

\begin{lemma}\label{New} Let $H$ be a finite group and let $x_1, y_1, x_2, y_2\in H$. Suppose that
\begin{enumerate}
\item[(1)] $o(x_1)$ and $o(y_1)$ are even;
\item[(2)] $\langle x_1^2, y_1^2, x_1y_1\rangle = H$;
\item[(3)] $\langle x_2, y_2\rangle = H$ and
\item[(4)] $o(x_1)o(y_1)o(x_1y_1)$ is coprime to $o(x_2)o(y_2)o(x_2y_2)$.
\end{enumerate}
If the above conditions are satisfied then $(G^0, x, y, g)$ is a mixed Beauville structure
for some $g\in(H\times H):(\mathbb{Z}/4\mathbb{Z})$ where $x = (x_1, x_2, 2), y = (y_1, y_2, 2)\in H\times H\times \mathbb{Z}/2\mathbb{Z}$ (note that $2\in\mathbb{Z}/4\mathbb{Z}$ generates the subgroup isomorphic to $\mathbb{Z}/2\mathbb{Z}$). Furthermore, if $H$ is a perfect group then we cane replace condition (2) with the condition
\begin{enumerate}
\item[(2')] $\langle x_1,y_1\rangle = H$.
\end{enumerate}
\end{lemma}

Note that in  \cite{8} this last hypothesis was incorrectly stated in terms of the perfectness of $G$ rather than $H$.
Bauer, Catanese and Grunewald go on to use the above lemma to construct
examples in the cases with the property that if $H$ is taken to be a sufficiently large alternating group
or a special linear groups $SL_2(p)$ with $p\not= 2, 3, 5, 17$ (though their argument also does not apply in the case $p=7$), then $(H\times H):(\mathbb{Z}/4\mathbb{Z})$ has
a mixed Beauville structure. Given the extent to which mixed Beauville groups are
in short supply it would be interesting to see if the above construction can be used
in other cases.

\begin{prob}Find other groups $H$ that the above lemma can be applied to.
\end{prob}

In \cite{NewP} the author and Pierro prove a slight generalisation of Lemma \ref{New} that replaces the cyclic group of order $4$ with the dicyclic group of order $4k$ defined by the presentation
\[
\langle x,y\,|\,x^{2k}=y^4=1, x^y=x^{-1}, x^k=y^2\rangle
\]
for some positive integer $k$. In particular, when finding examples of groups that satisfy the hypotheses of this generalisation (which is sufficient to show that such groups satisfy the hypotheses of Lemma \ref{New}) we obtain new examples of mixed Beuville groups from the groups $H$ and $H\times H$ where $H$ is any of the alternating groups $A_n$ ($n\geq6$), the linear groups $PSL_2(q)$ ($q\geq7$ odd), the unitary groups $PSU_3(q)$ ($q\geq3$), the Suzuki groups $^2B_2(2^{2n+1})$ ($n\geq1$), the small Ree groups $^2G_2(3^{2n+1})$ ($n\geq1$), the large Ree groups $^2F_4(q)$ ($q\geq8$), the Steinberg triality groups $^3D_4(q)$ ($q\geq2$) and the sporadic simple groups (including the Tits group $^2$F$_4(2)'$) as well as the groups $PSL_2(2^n)\times PSL_2(2^n)$ ($n\geq3$).

What about $p$-groups? If $p$ is odd then again, the absence of index 2 subgroups
ensures that there exist no mixed Beauville $p$-groups. In the construction described above the technical
constraints on $H$ ensure that it cannot be a 2-group, stopping this providing a
source of examples. Early examples of mixed Beauville 2-groups were given by Bauer, Cataneses and Grunewald constructed in \cite{10} where they constructed two mixed Beauville groups of order $2^8$. Even so, the lack of known  Beauville 2-groups makes the following a natural problem.

\begin{prob} Construct infinitely many mixed Beauville 2-groups.
\end{prob}

\section{Miscellanea}

\subsection{$PSL_2(q)$ and $PGL_2(q)$}\label{GarSec}

In \cite[Question 7.7]{9} Bauer, Catanese and Grunewald asked the following question,

\begin{quotation}
Existence and classification of Beauville surfaces, i.e.,

a) which finite groups $G$ can occur?

b) classify all possible Beauville surfaces for a given finite group $G$.
\end{quotation}

In \cite{28} Garion answered the above in the case of the groups $PSL_2(q)$ and
$PGL_2(q)$. For $PSL_2(q)$ we have the following.

\begin{thm}
Let $G = PSL_2(q)$ where $5 < q = p^e$ for some prime number $p$ and
some positive integer $e$. Let $\tau_1 = (r_1, s_1, t_1)$, $\tau_2 = (r_2, s_2, t_2)$ be two hyperbolic
triples of integers. Then $G$ admits an unmixed Beauville structure of type $(\tau_1, \tau_2)$
if, and only if, the following hold:
\begin{enumerate}
\item[(i)] the group $G$ is a quotient of the triangle groups $T_{r_1,s_1,t_1}$ and $T_{r_2,s_2,t_2}$ with
torsion-free kernel;
\item[(ii)] if $p = 2$ or $e$ is odd or $q = 9$, then $r_1s_1t_1$ is coprime to $r_2s_2t_2$. If $p$ is odd,
$e$ is even and $q > 9$, then $g = gcd(r_1s_1t_1, r_2s_2t_2)\in\{1,p,p^2\}$. Moreover, if $p$
divides $g$ and $\tau_1$ (respectively $\tau_2$) is up to a permutation $(p, p, n)$ then $n\not=p$
and $n$ is a `good $G$-order'.
\end{enumerate}
\end{thm}

Here by `good $G$-order' we mean the following. Let $q$ be an odd prime power
and let $n > 1$ be an integer. Then $n$ is a good $G$-order if either
\begin{itemize}
\item $n$ divides $(q-1)/2$ and a primitive root of unity $a$ of order $2n$ in $\mathbb{F}_q$ has the
property that $-a = c^2$ for some $c\in\mathbb{F}_q$ or
\item $n$ divides $(q + 1)/2$ and a primitive root of unity $a$ of order $2n$ in $\mathbb{F}_q^2$ has the
property that $-a = c^2$ for some $c\in\mathbb{F}_{q^2}$ such that $c^{q+1}=1$.
\end{itemize}

A similar theorem is given for the groups $PGL_2(q)$.

Given that generic lists of maximal subgroups of other low rank groups of Lie
type are well known in numerous other cases, it seems likely that analogous results
for these groups can also be obtained. We thus reiterate Bauer, Catanese and
Grunewald's earlier question in this case.

\begin{prob} Obtain results analogous to the above for other classes of finite
simple groups.
\end{prob}

\subsection{Fundamental Groups of Beauville Surfaces}\label{F1}

We mentioned in the introduction that Beauville surfaces have
fundamental groups that are easy to work with. To make this vague remark a little more specific we note the
following. Suppose that if $G$ is a Beauville group with a Beauville structure of type
$((a_1, b_1, c_1), (a_2, b_2, c_2))$, then for $i = 1, 2$ there exist surjective homomorphisms
$\rho_i:T_{a_i,b_i,c_i}\rightarrow G$. The direct product ker$(\rho_1)\times$ker$(\rho_2)$ is the fundamental group of the product $\mathcal{C}_1\times\mathcal{C}_2$. The fundamental group of the surface $(\mathcal{C}_1\times\mathcal{C}_2)/G$ is now an extension of  a normal subgroup ker$(\rho_1)\times$ker$(\rho_2)$ by $G$, or more precisely the inverse image in $T_{a_1,b_1,c_1}\times T_{a_2,b_2,c_2}$ of the diagonal subgroup of $G\times G$ under the epimorphism $\rho_1\times\rho_2$. It turns out that this simple description of the fundamental
group is responsible for the rigidity of Beauville surfaces and this in turn ensures that the
topological and geometric features of the surfaces are closely intertwined - see \cite[Section 9]{42} for details.

Unsurprisingly, since a Beauville group dictates so many features of its corresponding Beauville surface which
in turn determines its fundamental group we also have the reverse relationship
whereby the fundamental group determines the original Beauville group. The following is proved by Gonz\'{a}lez-Diez and Torres-Teigell \cite{32,48}. (It also worth noting related results given by Bauer, Catanese and Grunewald in \cite{9} and by Catanese in \cite{16}).

\begin{thm} Two Beauville surfaces are isometric if and only if their fundamental groups are isomorphic.
\end{thm}

The fundamental group is one of the most basic tools in algebraic topology. It is, however, somewhat limited in its usefulness and topologists have found several important higher dimensional analogues of the fundamental group and so it is natural to pose the following question.


\begin{qu}
Do the higher homotopy/homology/cohomology groups of a Beauville
surface have similar descriptions in terms of triangle groups and the corresponding Beauville group and to what extent do
they uniquely determine the surface?
\end{qu}

By way of partial progress on this question in \cite{BCF} Bauer, Catanese and Frapporti recently showed that for any Beauville surface $\mathcal{S}$ the homology group $H_1(\mathcal{S},\mathbb{Z})$ is finite. They also give a much more detailed discussion of geometric aspects of the study of fundamental groups of Beauville surfaces and related objects as well as computer calculations of these objects in some cases.

\subsection{Automorphism Groups of Beauville Surfaces}\label{F2}

In \cite{41} Jones investigated the automorphism groups of unmixed Beauville surfaces.
Some of these results were obtained independently by Fuertes and Gonz\'{a}lez-Diez
in \cite{23} and were later extended to mixed Beauville surfaces by Gonz\'{a}lez-Diez and
Torres-Teigell in \cite[Section 5.3]{31}.

\begin{thm}
The automorphism group Aut($\mathcal{S}$) of a Beauville surface $\mathcal{S}=(\mathcal{C}_1\times\mathcal{C}_2)/G$ has a normal subgroup Inn($\mathcal{S})\triangleleft Z(G)$ with Aut($\mathcal{S})/\mbox{Inn}(\mathcal{S})$ isomorphic to a subgroup of the wreath product $S_3\wr S_2$. In particular Aut($\mathcal{S}$) is a finite soluble
group of order dividing $72|Z(G)|$ and of derived length at most $4$.
\end{thm}

Here the subgroup Inn($\mathcal{S}$) consists of automorphisms preserving the two curves (or more precisely, induced by automorphisms of $\mathcal{C}\times\mathcal{C}'$ preserving them) though it does not necessarily contain all of them: they form a subgroup of index at most $2$ in Aut($\mathcal{S}$), whereas Inn($\mathcal{S}$) can have index up to 72. The results in the mixed case are similar.

\subsection{Beauville Genus Spectra}

In \cite[Question 7.7(b)]{9} Bauer, Catanese and Grunewald ask us to classify all
possible Beauville surfaces for a given finite group $G$.

As a partial answer to this, in \cite[Section 4]{25} Fuertes, Gonz\'{a}lez-Diez and Jaikin-Zapirain introduce the concept of Beauville genus spectrum which we define as
follows.

\begin{de} Let $G$ be a finite group. The \textbf{Beauville genus spectrum} of $G$
is the set $Spec(G)$ of pairs of integers $(g_1, g_2)$ such that $g_1\leq g_2$ and there are
curves $\mathcal{C}_1$ and $\mathcal{C}_2$ of genera $g_1$ and $g_2$ with an action of $G$ on $\mathcal{C}_1\times\mathcal{C}_2$ such that
$(\mathcal{C}_1\times\mathcal{C}_2)/G$ is a Beauville surface.
\end{de}

By the Riemann-Hurwitz formula each $g_i$ is bounded above by $1+\frac{|G|}{2}$
and so this set is always finite. Fuertes, Gonz\'{a}lez-Diez and Jaikin-Zapirain determine the
Beauville spectra of several small groups.

\begin{prop}
\begin{enumerate}
\item $Spec(S_5) = \{(19, 21)\}$
\item $Spec(PSL_2(7)) = \{(8, 49), (15, 49), (17, 22), (22, 33), (22, 49)\}$
\item $Spec(S_6) = \{(49, 91), (91, 121), (91, 169), (121, 169), (151, 169)\}$
\item If $gcd(n, 6) = 1$ and $n>1$ then $$Spec((\mathbb{Z}/n\mathbb{Z})\times(\mathbb{Z}/n\mathbb{Z})) = \{(\frac{(n-1)(n-2)}{2},\frac{(n-1)(n-2)}{2})\}.$$
\end{enumerate}
\end{prop}

Unpublished calculations of the author's PhD student, Emilio Pierro, has added a
few more finite simple and almost simple groups to the above list, the largest being the Mathieu group M$_{23}$. Furthermore, the Beauville genus structures of
$PSL_2(q)$ and $PGL_2(q)$ may be deduced from the results discussed in Subsection
\ref{GarSec}. This naturally leads us to ask the following.

\begin{prob} Determine the Beauville genus spectrum of more groups.
\end{prob}

\subsection{Characteristically Simple Groups}

Characteristically simple groups are usually defined in terms of characteristic subgroups, but for finite groups this turns out to be equivalent to the following.

\begin{de} A finite group $G$ is said to be \textbf{characteristically simple} if $G$ is isomorphic
to the direct product $H^k$ where $H$ is a finite simple group for some positive integer $k$.
\end{de}

If we fix $H$ then for large values of $k$ the group $H^k$ will not be 2-generated and
therefore will not be Beauville. For more modest values of $k$ there is, however, still
hope. These groups have recently been investigated by Jones in \cite{39, 40} where the
following conjecture is investigated.

\begin{conj} Let $G$ be a finite characteristically simple group. Then $G$ is
Beauville if and only if it is 2-generated and not isomorphic to the alternating
group A$_5$.
\end{conj}

Theorem \ref{SimpThm} shows that this conjecture is true for the characteristically simple
group $H^k$ in the case $k = 1$ for every non-abelian finite simple group $H$. If $G$ is abelian then this conjecture holds by Theorem \ref{AbThm} following the convention that a cyclic group is not considered to be 2-generated. In [39]
the above conjecture is verified for the alternating groups and in \cite{40} it is verified
for the linear groups $PSL_2(q)$ and $PSL_3(q)$, the unitary groups $PSU_3(q)$, the Suzuki groups
$^2B_2(2^{2n+1})$, the small Ree groups $^2G_2(3^{2n+1})$ and the sporadic simple groups. In
addition to the above the author has performed computations that verify the above
conjecture for all characteristically simple groups of order at most $10^{30}$. As an
amusing aside we note that this shows that whilst A$_5$ is not a Beauville group, the
direct product of nineteen copies of A$_5$ is!

In \cite{19} the author considers which of the characteristically simple groups are
strongly real Beauville groups. The main conjecture is the following.

\begin{conj} If $G$ is a finite simple group of order greater than 3, then $G\times G$ is a
strongly real Beauville group.
\end{conj}

It is likely that many larger direct products are also strongly real, however
the precise statement of a conjecture along these lines is likely to be much more
complicated. For example, a straightforward computation verifies that neither of
the groups M$_{11}\times$M$_{11}\times$M$_{11}$ and M$_{23}\times$M$_{23}\times$M$_{23}$ are strongly real despite the
fact that both of the groups M$_{11}\times$M$_{11}\times$M$_{11}\times$M$_{11}$ and M$_{23}\times$M$_{23}\times$M$_{23}\times$M$_{23}$
are.

The above conjecture has been verified for the alternating groups (though slightly stronger results are true in this case), the sporadic
simple groups, the linear groups $PSL_2(q)$ ($q > 5$), the Suzuki groups $^2B_2(2^{2n+1})$, the sporadic
simple groups (including the Mathieu groups $M_{11}$ and $M_{23}$, despite the statement of Conjecture 5.5) and all of the finite simple groups
of order at most $100\, 000\, 000$.

\subsection{Orbits of the Absolute Galois Group}

The task of understanding the absolute Galois group Gal($\overline{\mathbb{Q}}/\mathbb{Q}$) is of central importance in algebraic number theory and is related to the Inverse Galois Problem (it is
equivalent to asking if every finite group is a quotient of Gal($\overline{\mathbb{Q}}/\mathbb{Q}$) under a topologically closed normal subgroup) and this is arguably the hardest open problem
in algebra today. As things stand Gal($\overline{\mathbb{Q}}/\mathbb{Q}$) remains very poorly understood.
A natural approach to understanding any group is to study some action(s) of the
group. An immeidate consequence of Bely\u{\i}'s Theorem is that Gal($\overline{\mathbb{Q}}/\mathbb{Q}$) acts on the set of all Beauville surfaces.
Recently there has been much interest in constructing orbits consisting of mutually
non-homeomorphic pairs of Beauville surfaces. In \cite{32, 34} Gonz\'{a}lez-Diez, Jones and
Torres-Teigell have constructed arbitrarily large orbits of Gal($\overline{\mathbb{Q}}/\mathbb{Q}$) consisting of
mutually non-homeomorphic pairs of Beauville surfaces defined by the Beauville
groups $PSL_2(q)$ and $PGL_2(q)$.

\begin{prob} Construct arbitrarily large orbits of Gal($\overline{\mathbb{Q}}/\mathbb{Q}$) consisting of mutually non-homeomorphic pairs of Beauville surfaces using other groups.
\end{prob}

A slightly different motivation for addressing the above problem comes from
the following. Knowing whether or not Gal($\overline{\mathbb{Q}}/\mathbb{Q}$) acts faithfully on the set of
Beauville surfaces is equivalent to the longstanding question of whether or
not Gal($\overline{\mathbb{Q}}/\mathbb{Q}$) acts faithfully on the set of regular dessins. This was recently resolved by Gonz\'{a}lez-Diez and Jaikin-Zapirain in \cite{GGJZ} by showing that Gal($\overline{\mathbb{Q}}/\mathbb{Q}$) acts faithfully on the set of Beauville surfaces.

\end{document}